\documentclass[11pt]{amsart}
\usepackage{amsmath}
\usepackage{amsthm}
\usepackage{amssymb}
\usepackage{amscd}
\usepackage{amsfonts}
\usepackage{amsbsy}
\usepackage{epsfig,afterpage}
\usepackage[dvips]{psfrag}
\usepackage{subfigure}
\usepackage{epsfig}
\usepackage{amsmath}
\usepackage{srcltx}
\usepackage{amsfonts}
\usepackage{amssymb}
\usepackage{psfrag}
\usepackage{verbatim}
\usepackage{graphicx}
\usepackage{lscape}
\usepackage{scalefnt}
\usepackage{color,colortbl,multirow}
\usepackage{hhline}
\usepackage{overpic}
\usepackage[all]{xy}

\newcommand{\sgn}{\operatorname{sgn}}
\def\d{\delta}

\def\eps{\varepsilon}

\def\d{\delta}
\def\e{\varepsilon}

\newcommand{\bbox}{\ \hfill\rule[-1mm]{2mm}{3.2mm}}

\newtheorem {theorem} {Theorem}
\newtheorem {proposition}{Proposition}

\newtheorem {lemma}[theorem]  {Lemma}
\newtheorem {example} {Example}

\newcommand{\R}{\mathbb{R}}

\begin{document}

\title[Non--smooth Slow--Fast Systems]
{On Non--smooth Slow--Fast Systems}

\author[Jaime R. de Moraes and Paulo R. da Silva ]
{Jaime R. de Moraes $^1$ and Paulo R. da Silva $^2$}

\address{$^{1}$ Curso de Matem\'{a}tica -- UEMS,
Rodovia Dourados--Itaum Km 12, CEP 79804--970 Dourados, Mato Grosso do Sul,
Brazil.}

\email{jaime@uems.br}

\address{$^{2}$ Departamento de Matem\'{a}tica -- Instituto de Bioci\^{e}ncias Letras
e Ci\^{e}ncias Exatas, UNESP -- Univ Estadual Paulista, Rua C. Colombo, 2265,
CEP 15054--000 S\~{a}o Jos\'{e} do Rio Preto, S\~{a}o Paulo, Brazil}

\email{paulo.r.silva@unesp.br}


\date{}
\dedicatory{}

\maketitle
\begin{abstract} We deal with non--smooth differential systems 
	$\dot{z}=X(z), z\in\R^{n},$	with discontinuity occurring in 
	a codimension one smooth surface $\Sigma$. 
	A regularization  of $X$  is a 1--parameter family  of 
	smooth vector fields $X^{\delta},\delta>0$, satisfying that 
	$X^{\delta}$ converges pointwise to $X$ in 
	$\R^{n}\setminus\Sigma$,  when $\delta\rightarrow 0$. We work with two
	known regularizations: the classical one proposed by Sotomayor and Teixeira and 
	its generalization, using non-monotonic transition functions. 
	Using the techniques of geometric singular 
	perturbation theory we study  minimal sets of  regularized systems. 
	Moreover, non-smooth slow--fast systems are studied 
	and the persistence of the sliding region by singular perturbations is analyzed. 
\end{abstract}

\section{Introduction}
One finds in real life and in various branches of science distinguished phenomena 
whose mathematical models are expressed by  piecewise smooth systems and deserve 
a systematic analysis, see for instance \cite{Gua,Kuz,Med}. However sometimes the 
treatment of such objects is far from the usual techniques or 
methodologies found in the smooth universe. In fact, for such systems, everything 
we know from the qualitative theory of dynamical systems has its own 
versions, starting with the concept of solution.

\smallskip

Consider two smooth vector fields $X^+,X^-$ defined in $\R^n.$
A piecewise-smooth system  is $\dot{x}=X(x)$ with
\begin{equation}
\label{sis1} X=\dfrac{1}{2}\left[(1+\sgn (h))X^++(1-\sgn (h))X^-)\right],
\end{equation}
$h:\R^n\longrightarrow \R$ smooth and $0$ a regular value of $h$.
The set $\Sigma=\{x\in\R^n:h(x)=0\}$ is called 
\emph{switching manifold}.

\smallskip

In order to define what a solution is, it is necessary, first of all, to agree on 
what happens in $\Sigma$.
The points in  $\Sigma$ are classified  as \emph{regular} (if 
$X^+$ and $X^-$ are transversal to $\Sigma$) or \emph{singular} (if 
$X^+$ or $X^-$ is tangent to $\Sigma$). Moreover the regular points 
are classified according to Filippov's terminology \cite{AF} 
(\footnote{As usual, we denote  $Xf=\nabla f.X$.})
\begin{itemize}
	\item[(i)] $\Sigma^w =\{x\in\Sigma:(X^+h.X^-h)(x) >0 \}$ is the \textit{sewing} region;
	\smallskip
	\item[(ii)] $\Sigma^s=\{x\in\Sigma:(X^+h.X^-)h(x) <0 \}$ is the \textit{sliding} region.
\end{itemize}

To be more precise we subdivide $\Sigma^{s}$ in \textit{attracting sliding} 
$\Sigma^{s}_a$ (if $X^+h<0$ and $X^-h>0$) and \textit{repelling sliding or escape} 
$\Sigma^{s}_r$ (if $X^+h>0$ and $X^-h<0$).

\smallskip

The orbits of $X$ by $\Sigma^w$ are naturally concatenated. 
On $\Sigma^s$ is defined the \emph{sliding vector field} 
$X^{\Sigma}$ as a  linear convex combination of 
$X^+$ and $X^-$ which is tangent to $\Sigma$. The orbits by $\Sigma^s$ follow the
flow of $X^{\Sigma}$, a linear convex combination of $X^+$ and $X^-$
tangent to $\Sigma$, that is
\begin{equation}\label{GeralSVF}
X^{\Sigma} = \dfrac{( X^+ .h)X^--(X^- .h)X^+ }{(X^+-X^-).h}.
\end{equation}
The vector field $X^{\Sigma}$ is called \emph{sliding vector field}.

\smallskip

While Filippov said how the flow of a piecewise smooth vector field behaves 
when finding the set of discontinuity, Sotomayor and 
Teixeira (see \cite{ST}) addressed the problem by seeking smooth approximations 
which were called regularization.  A \textit{regularization}
of $X$ is a family of smooth vector 
fields $ X^{\delta} $ depending on a parameter $\delta > 0$ 
such that $ X^{\delta} $ converges simply to  $X$ in  $\R^{n}\setminus \Sigma $
when $\delta$ goes to zero.  

\smallskip

The Sotomayor-Teixeira regularization (\emph{ST-regularization})   is
the one parameter family $X^{\delta}$ given by 
\begin{equation}  \label{regst}
X^{\delta}= \Big(\frac{1+\varphi(h/\delta)}{2}\Big) X^+ + \Big(\frac{1-\varphi(h/\delta)}{2}\Big)X^-
\end{equation} 
where $\varphi:\R\rightarrow[-1,1]$ is a smooth function satisfying that 
$\varphi(t)=1$ for $t\geq 1$, $\varphi(t)=-1$ for $t\leq -1$ and $\varphi'(t)>0$ 
for $t\in(-1,1)$.
The regularization  is smooth for $\delta > 0$ and satisfies that 
$X^{\delta}=  X^+$ on $\{h\geq\delta\}$ and $X^{\delta}=  X^-$ on $\{h\leq-\delta\}$. 
The flow of the regularized vector field proposed by Sotomayor--Teixeira, after the limit 
process, is exactly the flow idealized by Filippov.

\smallskip

In 2005 Silva, Teixeira and Buzzi, strongly inspired by Freddy Dumortier, 
wrote the article \cite{BPT}. 
They proved that the regularization proposed by Sotomayor and Teixeira 
generates a singular perturbation 
problem. This provides a very important application of GSP--theory 
(geometric singular perturbation theory). 
Joint with Llibre they published  \cite{LST, LST3, LST4, LST5}. They 
studied regularization 
problems in $\R^n$, the double 
regularization in the case in which the discontinuity has codimension 1 
(intersection of two planes) and the regularization  in more degenerate 
surfaces. Bonet-Rev\'{e}s, Larrosa, M-Seara, Kristiansen,  
Uldall and  Hogan also studied the singular perturbation problem arising 
from regularization. More precisely, they analyzed the regularization of 
fold-fold singularities where bifurcation and canard boundary cycles can 
occur. In \cite {BRT, BLT}  the authors use of asymptotic analysis 
and the extension of the critical manifold to non--normally hyperbolic points. 
In \cite{Kr2} the authors strongly use the blow--up techniques developed by 
Dumortier--Roussarie \cite {DR} to deal with the same problem.

\smallskip

Our contribution to the general theory is to investigate the following problems.
\begin{itemize}
\item Applying techniques of the geometric singular perturbation theory we
study the  limit periodic sets (equilibrium points and 
periodic orbits  contained in the switching manifold) 
obtained as limit of orbits of  $X^{\d}$ when $\d \downarrow 0$  and we give 
alternative proofs of some results of \cite {SM}.  See \textbf{Section \ref{s2}}.
\smallskip
\item  The effect of breaking the monotonicity condition of the 
transition function $\varphi$ used in the regularization process \eqref{regst}.  
We review the concepts of sliding and sewing and their dependence on the regularization 
process considered.
What is the relation between such sets and those idealized by Filippov? 
See \textbf{Section \ref{s3}}.
\smallskip
\item Non--smooth slow--fast systems with sliding points in the critical manifold. 
We generalize the results of \cite{Ped1}
considering the new concept of sliding and sewing.
See  \textbf{Sections \ref{s4}} and  \textbf{\ref{s5}}.
\end{itemize}

In Section \ref{s2} we show how  the GSP-theory can be used to get information 
about the regularized vector field. From our previous papers we know that the 
trajectories of a  piecewise smooth vector field are obtained solving a 
slow--fast system with  critical manifold being the graphic of a smooth 
function defined in $\Sigma^{s}$. 
Moreover, the projection of the reduced flow, in $\Sigma$, is  the sliding flow. 
We prove that  hyperbolic equilibrium points $p$ (respectively periodic orbits 
$\gamma$) of the sliding vector field \eqref{GeralSVF} are limit of 
sequences of hyperbolic equilibrium points $p_{\delta}$ (respectively periodic 
orbits $\gamma_{\delta}$) of the regularized vector field \eqref{regst}. 
Besides, the dimensions of stable and unstable manifolds of $p_{\delta}$ 
(respectively $\gamma_{\delta}$) are determined. See \textbf{Theorem} \ref{teo3}.

\smallskip

\begin{figure}[h]\vspace{0.5cm}
\begin{overpic}[width=12cm]{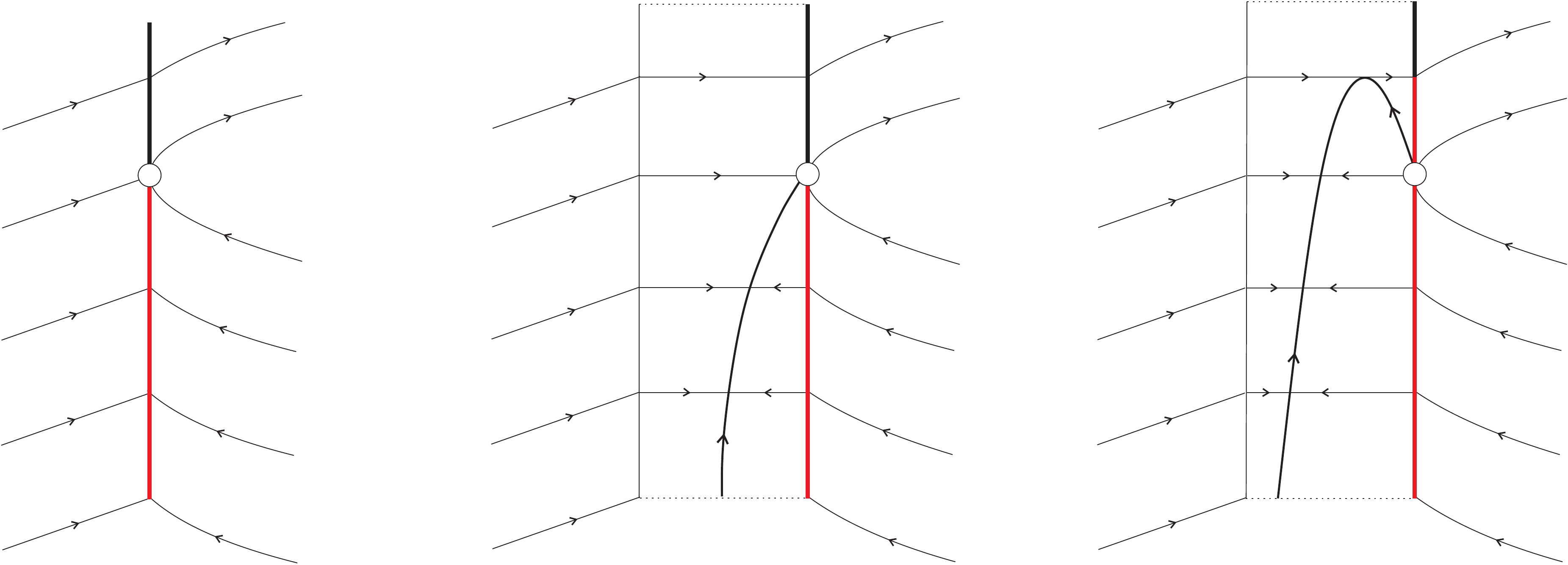}
		\put(7,37){(A)}
		\put(44,38){(B)}
		\put(83,38){(C)}
\end{overpic}
\caption{\footnotesize Figure (A) exhibits the flow of a piecewise smooth vector fields. 
Figure (B) is obtained using the ST--regularization. The sliding region is given 
in red. Figure (C) is obtained using a $r$--regularization. Note that 
$\Sigma^s\subseteq\Sigma_r^s$.}
\label{fig0}
\end{figure}

In Section \ref{s3} we propose a more general definition of sewing and 
sliding points (called $r$-\textit{sewing} and $r$-\textit{sliding} 
points, respectively). Sewing and sliding points  are defined depending on 
the choice of the regularization.  Roughly speaking, a point $p$ is a 
sewing point for a regularization $r$ if around $p$ the flow of $r$ 
is transversal to $\Sigma$. A point $p$ is a sliding point for $r$ if there exists 
a sequence of invariant manifolds of $r$ tending to a neighborhood of 
$p$ in $\Sigma$.  We prove that the sewing region $\Sigma^w$ contains the
	$r$-sewing region $\Sigma^w_r$ and the $r$-sliding region $\Sigma^s_r$ 
	contains the  sliding region $\Sigma^s$. Besides the sliding vector field 
	on $\Sigma^s$ can be smoothly extended on $\Sigma^s_r$.
	See Figure \ref{fig0} and  \textbf{Theorem} \ref{teoA}.
	
\smallskip

In Section \ref{s4} we study non--smooth slow--fast systems. 
Let $F,G:\R^{n+1}\times[0,+\infty)\rightarrow\R^{n}$, 
$H:\R^{n+1}\times[0,+\infty)\rightarrow\R$ and 
$h:\R^{n+1}\times[0,+\infty)\rightarrow\R$ be smooth.
We consider systems of the kind
\[\dot{x}=\left\{
\begin{array}{ll}
F(x,y,\eps), \quad\hbox{if} \quad h(x,y,\eps)\geq 0\\
G(x,y,\eps),\quad \hbox{if} \quad h(x,y,\eps)\leq 0
\end{array}
\right. , \quad \eps\dot{y}=H(x,y,\eps).\]
We prove that the $r$-sliding region in the critical manifold $H(x,y,0)=0$ persists 
for small $\e>0$. See  \textbf{Theorem} \ref{teoB}.

\smallskip

In Section \ref{s5} we give one more definition of sewing and sliding 
points using any continuous combination (not necessary convex) of 
$X^+$ and $X^-$. We state and prove that sliding regions obtained via continuous 
combination of $X^+$ and $X^-$ also are persistent by singular perturbation. 
See \textbf{Theorem} \ref{teoC}.

\section{Fenichel's Theory and Sliding Vector Fields.}\label{s2}
In order to simplify our explanation  we take local local 
coordinates such that
\begin{equation*}\label{localcoordinate}
\R^n=\R^{n-1}\times\R,\quad\Sigma=\{(x,y)\in\R^{n-1}\times\R:y=0\}. 
\end{equation*}
Let $X^+,X^-:\R^n\rightarrow\R^n$ be smooth vector fields. Denote  $X^+=(f_1,g_1)$ 
and $X^-=(f_2,g_2)$. The piecewise smooth vector field which we consider in 
this section is
\begin{equation}
\label{psvf} X=\dfrac{1}{2}\left[(1+\sgn (y))X^++(1-\sgn (y))X^-)\right].
\end{equation}

The sliding vector field \eqref{GeralSVF} becomes
\begin{equation}\label{svfsimples}
X^{\Sigma}=\dfrac{( X^+ .y)X^--(X^- .y)X^+ }{(X^+-X^-).y}
=\left(\dfrac{f_2g_1-f_1g_2}{g_1-g_2},0\right).
\end{equation}
The trajectories of  the ST-regularized vector field $X^\d$
given by \eqref{regst} are the the solutions of the differential system
\begin{equation}\label{sisreg}
\mathcal{P}_{\delta}:\left\{
\begin{array}{lll}
\dot{x}&=&(f_1+f_2)/2+\varphi(y/\delta)(f_1-f_2)/2\\
\dot{y}&=&(g_1+g_2)/2+\varphi(y/\delta)(g_1-g_2)/2
\end{array}
\right..
\end{equation}

Consider the polar blow up  $(x,y,\delta)=\Phi(x,\theta,r)= (x,r\cos\theta,r\sin\theta)$.
We get the system
\begin{equation*}\label{polar}
\bar{\mathcal{P}_{r}}:
\left\{
\begin{array}{lcl}
\,\,\,\dot x &=&\left(f_1+f_2\right) /2 +\varphi\left
(\cot \theta \right)  \left( f_1-f_2\right) /2, \\
r\dot{\theta}&=&-\sin \theta\left[ \left( g_1+g_2\right)
/2 +\varphi\left(\cot \theta \right)\left( g_1-g_2\right) /2\right].
\end{array}
\right.
\end{equation*} 
For $r=0$, $\bar{\mathcal{P}_{r}}$ has two limit problems, the reduced \eqref{polar slow}  
and the layer, as we will see below.
\begin{equation}\label{polar slow}
\bar{\mathcal{P}_{0}}^{reduced}:
\left\{
\begin{array}{rcl}
\,\,\,\dot x &=&\left(f_1+f_2\right) /2 +\varphi\left
(\cot \theta \right)  \left( f_1-f_2\right) /2, \\
0&=&-\sin \theta\left[ \left( g_1+g_2\right)
/2 +\varphi\left(\cot \theta \right)\left( g_1-g_2\right) /2\right].
\end{array}
\right.
\end{equation} 
\begin{equation*}\label{polar fast}
\bar{\mathcal{P}_{0}}^{layer}:
\left\{
\begin{array}{rcl}
\,\,\,x '&=&0, \\
\theta'&=&-\sin \theta\left[ \left( g_1+g_2\right)
/2 +\varphi\left(\cot \theta \right)\left( g_1-g_2\right) /2\right].
\end{array}
\right.
\end{equation*}

The polar blow--up is not  the most suitable for calculations.
For this reason, we perform the directional blow-up. The directional blow--up consists in the
following change of coordinates  $(x,y,\delta)=\Gamma(x,\overline{y},\delta)=
(x,\delta\overline{y},\delta).$ 
We observe that the direcional blow--up and the polar blow--up 
are essentially the same. In fact, if we consider the map
$G(x,\theta,r)=(x,\cot\theta, r\sin\theta)$ then $\Gamma\circ G=\Phi.$
The direcional blow--up applied in  system \eqref{sisreg} gives
\begin{equation}
\label{regblow}
\bar{\mathcal{P}_{\delta}}:\quad
\dot{x}=\alpha(x,\bar{y},\delta),\quad
\delta\dot{\bar{y}}=\beta(x,\bar{y},\delta)
\end{equation}
with  \[\alpha=(f_1+f_2)/2+\varphi(\bar{y})(f_1-f_2)/2,\quad
\beta=(g_1+g_2)/2+\varphi(\bar{y})(g_1-g_2)/2\] and the $f_1,f_2,g_1,g_2$  evaluated 
at $(x,\delta\bar{y})$.

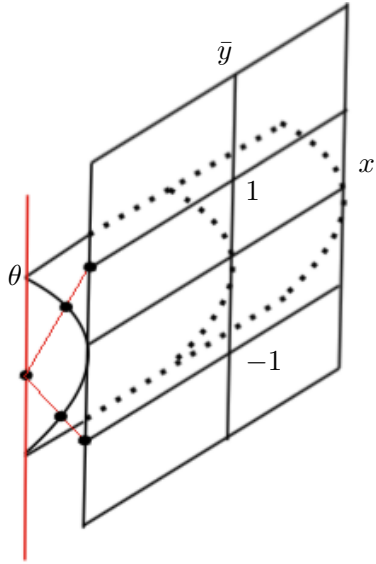
\begin{figure}
	\begin{center}
		\begin{overpic}[width=6cm]{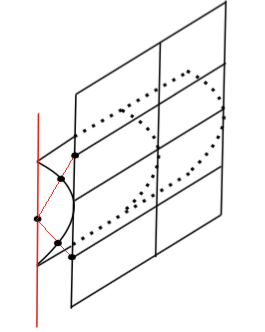}
			\put(50,35){$-1$}\put(50,65){$1$}\put(70,70){$x$}\put(45,90){$\bar{y}$}\put(8,50){$\theta$}
		\end{overpic}
		
	\end{center}
	\caption{\footnotesize Correspondence between the polar coordinates $(r, \theta) $ 
		and the directional coordinates $ (x, \bar {y}) $. The phase
		portrait on the semi--cylinder  $r=0,(\theta,x)\in[0,\pi]\times\R^{n-1}$ is 
		the central projection of the phase portrait of system \eqref{regblow} with
		$\delta=0,\bar{y}\in[-1,1]$.}
\end{figure}

\subsection{GSP-theory.} 
Systems as  \eqref{regblow} are known in the literature as  \emph{slow--fast systems}. 
General slow--fast systems are systems of the kind
\begin{equation}\label{pertub}
\dot{x} =\alpha(x,y,\d),\quad \d\dot{y}=\beta(x,y,\d),
\end{equation}
where $x\in\R^n$, $y\in\R^k$, $\d\geq0$ and $\alpha$ and $\beta$ are smooth functions.
Taking $\d=0$ in \eqref{pertub} we obtain the \emph{reduced system}
\begin{equation*}\label{reduzido}
\dot{x} =\alpha(x,y,0),\quad \beta(x,y,0)=0.
\end{equation*}
The set $\mathcal{S}=\{\beta(x,y,0)=0\}$ is called \emph{critical (or slow) 
manifold}. The time scale $\tau =t/\d$ transforms system \eqref{pertub} 
in the \emph{fast system}
\begin{equation}\label{fast}
\dot{x} =\d \alpha(x,y,\d),\quad \dot{y}=\beta(x,y,\d).
\end{equation}
Taking $\d =0$ in \eqref{fast} we get the \emph{layer system}.
We say that a point $(x_0,y_0)\in \mathcal{S}$ is \emph{normally hyperbolic}
if the real parts of the eigenvalues of $D_y\,\beta(x_0,y_0,0)$ are nonzero.

\smallskip

Let $\mathcal{N}\subset \mathcal{S}$ be a compact normally hyperbolic set. 
Consider $(x,y)\in\mathcal{N}$
and suppose that
$D_y\,\beta(x,y,0)$ has $k^s$ eigenvalues with negative real parts and
$k^u$ eigenvalues with positive real parts.
The following result ensures the persistence of normally
hyperbolic sets in $\mathcal{S}$ as invariant manifolds of system 
\eqref{pertub}, for small values of $\d>0$.

\begin{proposition}[Fenichel, \cite{F}]\label{Fenichel}
	Let $\mathcal{N}\subset \mathcal{S}$ be a $j$--dimensional compact
	normally hyperbolic manifold with a $(j+j^s)$--dimensional local stable 
	manifold $\mathcal{W}^s$ and a $(j+j^u)$--dimensional local unstable 
	manifold $\mathcal{W}^u$. Then there exists a family $\mathcal{N}_\d$ 
	such that the following statements hold.
	\begin{itemize}
		\item[(a)] $\mathcal{N}_0=\mathcal{N}$.
		
		\smallskip
		
		\item[(b)] $\mathcal{N}_\d$ is an invariant manifold of \eqref{pertub}
		with a $(j+j^s+k^s)$--dimensional local stable manifold $\mathcal{N}^s_\d$
		and a $(j+j^u+k^u)$--dimensional local unstable manifold $\mathcal{N}^u_\d$.
	\end{itemize}
\end{proposition}

\begin{proposition}\label{teo2} 
Let $X$ be the piecewise smooth vector field \eqref{psvf} and $X^{\d}$ its 
ST-regularization \eqref{regst}. If $\Sigma=\Sigma^{sl}$ then all points on 
the critical manifold $\mathcal{S}$ of system \eqref{regblow} are  normally 
hyperbolic. In particular $\mathcal{S}$ is a graphic of a smooth function 
$\bar{y}=h(x)$ with $(x,0)\in\Sigma^{sl}$. Moreover, the projection of the 
reduced flow,  on $\Sigma$, is  the sliding flow of \eqref{svfsimples}.
\end{proposition}

\noindent \emph{Proof.} For completeness of the text we rewrite the proof
originally presented in \cite{LST2}. The critical manifold $\mathcal {S}$ is 
defined by the equation $(g_1+g_2)+\varphi (\bar{y})(g_1- g_2)=0,$
with $g_1,g_2$ evaluated at $(x,0)$. If $(x,0) \in\Sigma^{s}$ then $g_1g_2<0$ and thus
\[(x,\bar{y})\in\mathcal{S}\iff \varphi (\bar{y})=-\dfrac{g_1+g_2}{g_1-g_2}.\] 
Since 
$|\frac{g_1+g_2}{g_1-g_2}|\leq 1$ and  $|\frac{g_1+g_2}{g_1-g_2}|=1$ if and 
only if  $g_1g_2=0$ it follows that 
$D_{\bar{y}} \beta(x,\bar{y},0)=\varphi'(\bar{y})(g_1-g_2)(x,0)\neq0.$
In fact,  for $(x,0) \in\Sigma^{s}$ we have $|\frac{g_1+g_2}{g_1-g_2}|<1$, 
$\bar{y}\in(-1,1)$ and thus $\varphi'(\bar{y})>0.$ It concludes the proof 
of the assertion about the critical manifold.

\smallskip

To see the relation between the reduced flow and the sliding vector vector 
field observe that taking $\varphi(\bar{y})=-\frac{g_1+g_2}{g_1-g_2}$ the 
system  \eqref{regblow} becomes 
\[\dot{x}=\frac{f_2g_1-f_1g_2}{g_1-g_2}\]  
which is exactly the same equation of the trajectories of $X^{\Sigma}$. \bbox

\begin{theorem}\label{teo3}  
Let $X$ be the piecewise smooth vector field \eqref{psvf}. If $\mathcal{Q}$ is a 
$\ell-$ dimensional compact invariant manifold of $X^{\Sigma}$  given by 
\eqref{svfsimples} with a $\ell+ \ell^s-$ dimensional local stable manifold, 
a $\ell+\ell^u-$ dimensional local unstable manifold and $\ell^s+\ell^u=n-1$. 
If $\mathcal{Q}\in\Sigma^{s}_a$  then there are a neighborhood $V$ of 
$\mathcal{Q}$ in $\R^{n}$ and $\d_0>0$ such that for $0 <\d<\d_0$, $ X^{\d}$ 
has a invariant manifold  $\mathcal{Q}_{\d}\in V$ with  
$(\ell+ \ell^s+1)-$dimensional stable manifold and $(\ell+\ell^u)-$dimensional 
unstable manifold.
\end{theorem}

\noindent\textit{Proof.} Assume $\ell=0$. Thus  $\mathcal{Q}=(x_0,0)$ is a 
hyperbolic equilibrium point of $X^{\Sigma}$. The trajectories of the 
ST-regularization $X^{\d}$ are the solutions of system \eqref{sisreg}. 
Performing the directional blow up we get 
system \eqref{regblow}. Moreover, $\mathcal{Q}$ is an equilibrium of 
$\dot{x}=\frac{f_2g_1-f_1g_2}{g_1-g_2}$, $\bar{y}=h(x)$
with $\mathcal{S} :\bar{y}=h(x)$ defined implictly by 
$(g_1+g_2)/2+\varphi(\bar{y})(g_1-g_2)/2=0.$
The hypothesis $g_1(\mathcal{Q})<0$ and $g_2(\mathcal{Q})>0$ implies that 
$\mathcal{S} $ satisfies the attractiveness condition:
\begin{itemize}
	\item all points in $\mathcal{S}$ are normally hyperbolic;
	\smallskip
	\item $k^s=1$ and $k^u=0$.
\end{itemize}
Since the reduced system and the sliding system have the same equations we 
conclude that $(x_0,y_0)$, with $y_0=h(x_0)$, is an equilibrium of the 
reduced system with $j^s=\ell^s$ and $j^u=\ell^u$. 
Then, applying Proposition \ref{Fenichel} we conclude that there exists 
$\mathcal{Q}_{\d}$, an equilibrium point of system \eqref{sisreg}, with  
$( \ell^s+ 1)$-dimensional stable 
manifold and $(\ell^u)$-dimensional unstable manifolds.  \bbox\\

\begin{example}\rm 
Let $X:\R^2\rightarrow\R^2$ be the piecewise smooth vector 
field \eqref{psvf} with  $X^+=(0,-1)$ and  $X^-=(x,-y+1)$.
The sliding vector field  \eqref{svfsimples} is 
\[X^{\Sigma}=\left(\dfrac{x}{2},0\right).\]
$\mathcal{Q}=(0,0)\in\Sigma^{s}_a$ is an repelling equilibrium point of 
$X^{\Sigma}$, that is,  $\ell^s=0$ and $\ell^u=1.$ 
Thus, Theorem \ref{teo3} says that for small $\delta$, $ X^{\d}$ has 
an equilibrium $\mathcal{Q}_{\d}$ of the kind saddle. 

\smallskip

Note that this can be verified directly with simple calculations. In fact, 
\[X^{\delta}=\left(\dfrac{x}{2}\left(1-\varphi\left(\dfrac{y}{\delta}\right)\right),
-\dfrac{y}{2}+ \varphi\left(\dfrac{y}{\delta}\right)\left(\dfrac{y-2}{2}\right)  \right)\]
and 
\[\mathcal{Q}_{\d}= (0,y_0), \quad \mbox{with}\quad \varphi\left(\dfrac{y_0}{\delta}\right)=\dfrac{y_0}{y_0-2}.\]
The existence of $\mathcal{Q}_{\d}$ is guaranteed by the fact that the 
graphics of $\varphi(y/\delta)$ and $\dfrac{y}{y-2}$
intersect at some value of $y\in (-1, 1)$. Moreover the eigenvalues of 
the linearized system at $(0,y_0)$ are 
\[\lambda_1=\dfrac{1}{2},\quad \lambda_2=
\dfrac{1}{y_0-2}+\varphi'\left(\dfrac{y_0}{\delta}\right).\dfrac{1}{\delta}\dfrac{y_0-2}{2}.\]
Since $y_0\in(-1,1)$ we have $\lambda_2<0$. Soon $\mathcal{Q}_{\d}$ is a
saddle. 
\end{example}

\begin{example}\rm 
Let $X:\R^3\rightarrow\R^3$ be the piecewise smooth 
vector field \eqref{psvf} with  
\[X^+=\left(0,2x_1+2x_2(\sqrt{x_1^2+x_2^2}-1),-1\right)\] and  
\[X^-=\left(-2x_2+2x_1(\sqrt{x_1^2+x_2^2}-1),0,1\right).\]
The sliding region is $\Sigma=(x_1,x_2)$--plane  and the sliding vector 
field is
\[X^{\Sigma}=\left(-x_2+x_1(\sqrt{x_1^2+x_2^2}-1) , 
x_1+x_2(\sqrt{x_1^2+x_2^2}-1),0\right).\]
The equilibrium point $p_0=(0,0,0)$ and the limit cycle 
$\gamma_0: x_1^2+x_2^2=1$, of $X^{\Sigma}$ satisfy the hypotheses of  
Theorem \ref{teo3} with $\ell=0$ and $\ell=1$ respectively.
\end{example}

\section{Non--smooth Systems and Regularization}\label{s3}
Consider piecewise smooth system \eqref{sis1} defined in an open set 
$\mathcal{U}\subset \R^n $ with
$X^{+},X^{-}:\mathcal{U}\rightarrow \R^n$, 
$h:\mathcal{U}\longrightarrow \R$ smooth and assume that $0$ is a regular value 
of $h$. 

\smallskip

The following regularization  will be refered as 
$r$--\textit{regularization}:
\begin{equation}\label{regtran}
X^{\delta}_r=\dfrac{\left(1+\psi\left(p,h/\delta\right)\right)}{2}X^+
+\dfrac{\left(1-\psi\left(p,h/\delta\right)\right)}{2}X^-,
\end{equation}
where $\psi:\Sigma \times\R\rightarrow [-1,1]$ is a more general smooth 
transition function satisfying that $\psi(p,t)=-1$ for $t\leq-1$ and 
$\psi(p,t)=1$ for $t\geq1$.

\smallskip

The definitions of $r$-sewing and $r$-sliding points were introduced 
in \cite{PR} and depend of the regularization $r$ considered. 
More precisely, $p\in\Sigma$ is a 
$r$-\emph{sewing} point if there exist an open neighborhood 
$\mathcal{U}\subset \R^{n-1}\times\R$ of $p$ and local coordinates $(x,y)$ 
defined in $\mathcal{U}$ such that:
\begin{itemize}
	\item[(a)] $\Sigma=\{y=0\}$;
	
	\smallskip
	
	\item[(b)] for each sufficiently small $\delta>0$, the 
	vector field $v(x,y)=(0,1)$  is a generator of \eqref{regtran} 
	in $\mathcal{U}$.
\end{itemize}
$p\in \Sigma$ is a $r$-\emph{sliding} point  if  there exist an open 
neighborhood $\mathcal{U}\subset \R^{n}$ of $p$ and a family of smooth 
manifolds $S_{\delta}\subset \mathcal{U}$ satisfying:
\begin{itemize}
	\item[(a)] $S_{\delta}$ is invariant for $r$;
	
	\smallskip
	
	\item[(b)] for each compact $K\subset\mathcal{U}$, the sequence 
	$S_{\delta}\cap K$ converges to $\Sigma \cap K$, as $\delta$ goes to 
	zero according to Hausdorff distance.
\end{itemize}

We denote $\Sigma^{w}_{r}$ and $\Sigma^{s}_{r}$ the $r$-sewing and 
$r$-sliding regions, respectively and
we assume local coordinates $(x,y)\in\R^{n-1}\times\R$ such that $h(x,y)=y$,
$X^+=(f_1,g_1)$,  $X^-=(f_2,g_2)$ with \[(f_i(x,y),g_i(x,y))\in\R^{n-1}\times\R,\quad i=1,2.\]
\begin{figure}[h]\vspace{0.3cm}
	\begin{overpic}[width=3cm]{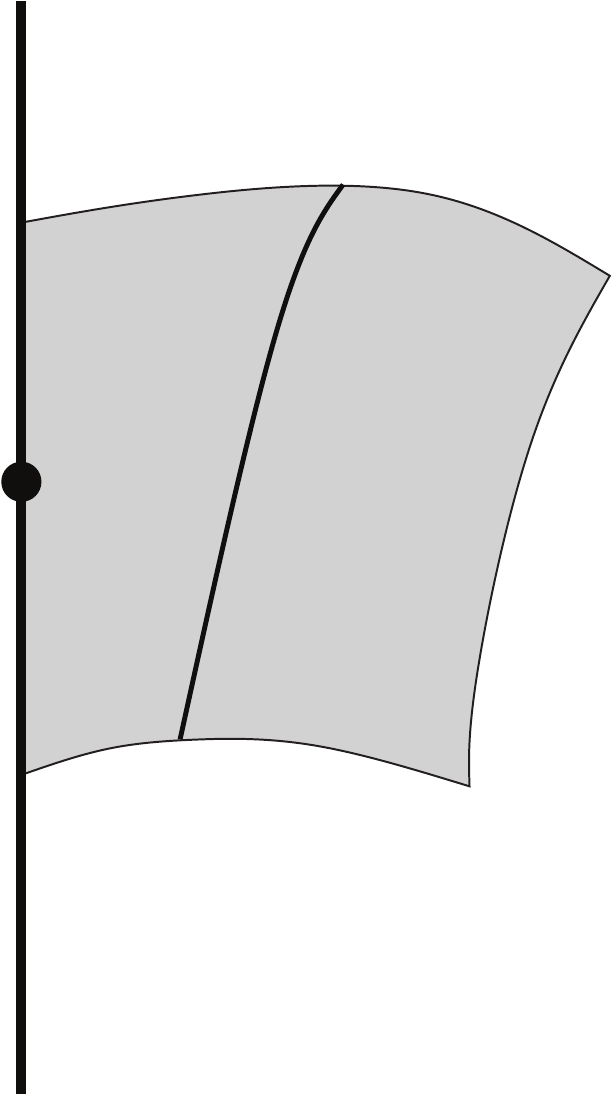}
		\put(-8,51){ $p$}
		\put(10,55){ $\mathcal{S}_\delta$}
		\put(-1,102){ $\Sigma$}
	\end{overpic}
	\caption{\footnotesize A $r$-sliding point $p$ defined by $r$-regularization.}
	\label{figinv}
\end{figure}

\begin{theorem}\label{teoA} Consider a non--smooth system  \eqref{sis1} and a 
	$r$-regularization \eqref{regtran}. Suppose that 
	$\frac{ \partial \psi }{\partial t}(x,t)\neq0$, for all $x\in\Sigma^s$ 
	and $|t|<1$. Then the following statements hold.
\begin{itemize}
	\item[(a)] $\Sigma_r^w\subseteq\Sigma^w$ and
	$\Sigma^s\subseteq \Sigma^s_r$.
	\smallskip
	\item[(b)] If $g_1\neq g_2$ in $\Sigma_r^s\setminus\Sigma^s$, 
		then the sliding vector field $X^\Sigma$ can be smoothly extended from $\Sigma^s$ to
		$\Sigma_r^s$.
\end{itemize}\end{theorem}
\noindent{\textit {Proof.}}
Consider $x_0\in\Sigma^{w}_{r}$. Since 
$v(x,y)=(0,1)$ is a generator of  \eqref{regtran}, $X^{-}(x_0,0)$ and $X^+(x_0,0)$ 
point to the same hand side. So \[(X^+. y)(X^-. y)(x_0,0)>0\] and we 
conclude that $x_0\in\Sigma^{w}$. Therefore 
$\Sigma^{w}_{r}\subseteq\Sigma^{w}$.

\smallskip

Now, consider $x_0\in\Sigma^{s}$. Taking the directional blow--up
$y=\delta \overline{y}$, the $r$-regularization \eqref{regtran} 
becomes the slow--fast system
\begin{equation}\label{eqth0}
\dot{x}=\alpha(x,\overline{y},\delta),\quad\delta\dot{\overline{y}}=\beta(x,\overline{y},\delta)
\end{equation}
where
\[
\alpha(x,\overline{y},\delta)=1/2\big((1
+\psi(x,\overline{y}))f_1+
(1-\psi(x,\overline{y}))\big)f_2
\]
and 
\[
\beta(x,\overline{y},\delta)=1/2\big((1
+\psi(x,\overline{y}))g_1+
(1-\psi(x,\overline{y}))\big)g_2,
\]
with the functions $f_i$ and $g_i$, $i=1,2$, evaluated 
at $(x,\delta \overline{y})$. The reduced system 
\begin{equation}\label{redblow}
\dot{x}=\alpha(x,\overline{y},0),
\end{equation}
is defined for $(x,\overline{y})$ in the critical manifold 
$\mathcal{S}=\{\beta(x,\overline{y},0)=0\}$.
Consider $|y_0|<1$ such that $(x_0,y_0)\in \mathcal{S}$. 
Note that $(g_1- g_2)(x_0,0)\neq0$ because $g_1g_2(x_0,0)<0$.  
So, 
\[\dfrac{\partial \beta}{\partial \overline{y}}(x_0, y_0,0)=
\dfrac{\partial \psi}{\partial \overline{y}}(x_0,y_0).
(g_1-g_2)(x_0,0)\neq0\]
and thus $(x_0,y_0)$ is a normally hyperbolic point. The Fenichel's
result ensures the existence of an invariant manifold
$\mathcal{S}_{\delta}$ of \eqref{eqth0} such that 
$\mathcal{S}_{\delta}\rightarrow \mathcal{S}_0$ as 
$\delta\rightarrow 0$, according to Hausdorff's distance. 
Thus $x_0\in\Sigma_r^s$ and it concludes the proof of item (a).

\smallskip

Since $(g_1 - g_2)(x,0)\neq0$ for  $x\in\Sigma^s_r$, the equation 
$\beta(x,\overline{y},0)=0$ provides
\begin{equation}\label{phi}
\psi\left(x,\overline{y}\right)=-\dfrac{g_1+g_2}{g_1-g_2}.
\end{equation}
The dynamics at $(x,\overline{y})$ obtained from \eqref{redblow} and 
\eqref{phi} is given by
\[\dot{x}=\dfrac{f_2g_1-f_1g_2}{g_1-g_2}.\]
Hence $X^\Sigma$ can be smoothly extended on 
$\Sigma_r^s$. It concludes the proof of item (b). \bbox\\

The following example illustrates Theorem \ref{teoA}.

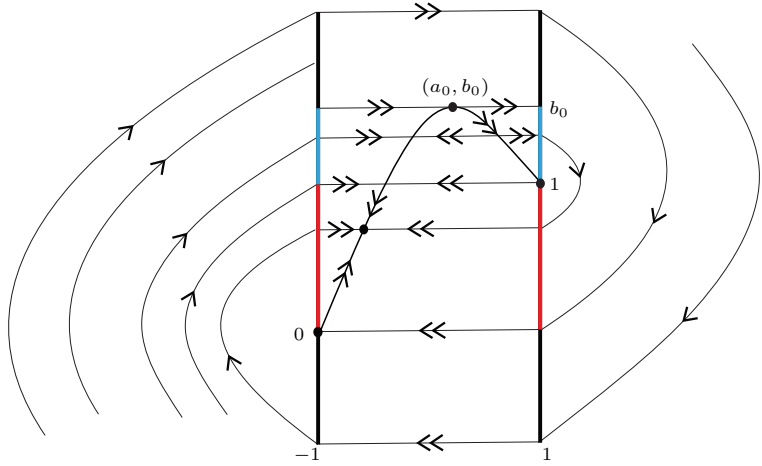
\begin{figure}[h]
	\begin{overpic}[width=10cm]{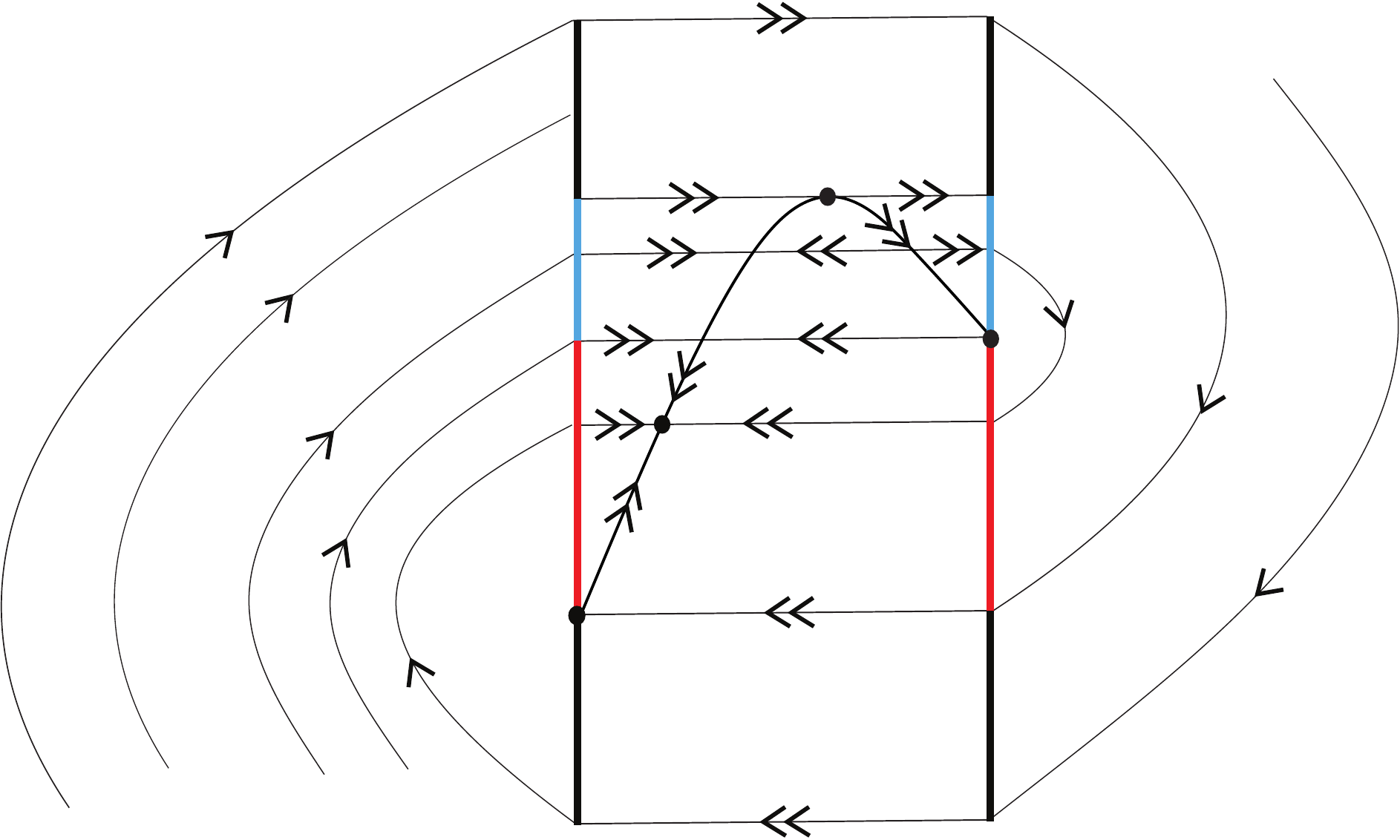}
		\put(38,-1){\tiny $-1$}
		\put(71,-1){\tiny $1$}
		\put(55,48){\tiny $(a_0,b_0)$}
		\put(38,15){\tiny$ 0$}
		\put(72,35){\tiny$ 1$}
		\put(72,45){\tiny$ b_0$}
	\end{overpic}
	\caption{\footnotesize Slow and fast dynamics of Example \eqref{exblow}.
		The black region is the $r$-sewing region, the red region is the
		sliding region according to Filippov convention. The sliding region in 
		blue appears considering the $r$-regularization  of system 
		\eqref{eqexblow}.}
	\label{fig00}	
\end{figure}

\begin{example}\rm \label{exblow}
	Consider the non--smooth system
	\begin{equation}\label{eqexblow}
	\dot{x}=(\dot{x}_1,\dot{x}_2)=\left\{
	\begin{array}{ccc}
	(x_2-1,-1),& \hbox{if}&  x_1\geq 0,\\
	(x_2,1),& \hbox{if} & x_1\leq 0.
	\end{array}
	\right.
	\end{equation}
	The sliding region is $\Sigma^s=]0,1[$ and the sliding vector field 
	\eqref{GeralSVF} is \[X^{\Sigma}=(0,-2x_2+1).\]
	The $r$-regularization of \eqref{eqexblow}
	is the 1--parameter family
	\[\dot{x_1}=1/2\big(2x_2-1-\psi\left(x_1/\delta,x_2\right)\big),
	\quad \dot{x_2}=-\psi\left(x_1/\delta,x_2\right).\]
	Suppose that partial derivative $\frac{\partial \psi}{\partial t}(t,x_2)$
	vanishes only at $(t,x_2)=(a_0,b_0)$ with $-1<a_0<1$ and $b_0>1$.
	Applying the directional blow--up $x_1=\delta\overline{x}_1$  we 
	obtain the slow--fast system
	\begin{equation}\label{regex1}
	\delta\dot{\overline{x}}_1=1/2\big(2x_2-1-\psi\left(
	\overline{x}_1,x_2\right)\big),\quad \dot{x}_2=-\psi(\overline{x}_1,x_2).
	\end{equation}
	The critical manifold $\mathcal{S}$ is given implicitly by
	\[x_2=\dfrac{\psi(\overline{x}_1,x_2)+1}{2}\]
	and it is a curve connecting $(\overline{x}_1,x_2)=(-1,0)$ to
	$(\overline{x}_1,x_2)=(1,1)$. All points in 
	$\mathcal{S}_0=\mathcal{S}\setminus{ (a_0,b_0)}$
	are normally hyperbolic.
	So Fenichel's result ensures the existence of an invariant
	manifold $\mathcal{S}_{\delta}$ of \eqref{regex1} converging to 
	$\mathcal{S}_0$ according to Hausdorff's distance.
	Thus $\Sigma_r^s=\left]0,b_0\right[$. The reduced system in $\mathcal{S}_0$
	is 
	\[\dot{x}_2=-2x_2+1.\]
	Thus $X^\Sigma$ can be smoothly extended on $\Sigma_r^s=\left] 0,b_0\right[  $.
	Note that the $r$-sliding region  $\Sigma_r^s$ obtained 
	contains the sliding region $\Sigma^s=\left]0,1\right[ $ defined by Filippov.
	See Figure \ref{fig00}.
\end{example}

\section{Non--smooth Slow--Fast Systems} \label{s4}
A non--smooth slow--fast 
system is 
\begin{equation}\label{eq1}
\dot{x}=\left\{
\begin{array}{ll}
F(x,y,\eps), \quad\hbox{if} \quad h(x,y,\eps)\geq 0\\
G(x,y,\eps),\quad \hbox{if} \quad h(x,y,\eps)\leq 0
\end{array}
\right. , \quad \eps\dot{y}=H(x,y,\eps)
\end{equation}
where $x\in\R^n$, $y\in\R$ and $\eps>0$ is a small parameter.
For each $\eps \geq 0$ let $\Sigma_{\eps}$ be the switching manifold
of \eqref{eq1}, i.e. $\Sigma_{\eps} = \{h(x,y,\eps) = 0\}$. 
We assume that $\Sigma_0$ and $\mathcal{S}=\{H(x,y,0) = 0\}$ are transversal and 
$H_y(x,y,0)\neq0$ on $\mathcal{S}$.
Taking $\e=0$ in \eqref{eq1} we have the \emph{reduced system}
\begin{equation}\label{eq22}
\dot{x} = \left\{
\begin{array}{l} \widetilde{F}(x) \quad $if$ \quad \tilde{h}(x) \geq 0, \\
\widetilde{G}(x)\quad $if$ \quad \tilde{h}(x) \leq 0,
\end{array}\right. \qquad  \widetilde{H}(x)=0,
\end{equation}
where $\widetilde{F}(x) = F(x,y(x),0)$, $\widetilde{G}(x)=G(x,y(x),0)$,
$\tilde{h}(x) = h(x,y(x),0)$ and $\tilde{H}(x) = H(x,y(x),0).$ System \eqref{eq22}  is 
defined in the critical manifold 
$\mathcal{S} = \{\widetilde{H}(x) = 0\}$.

\smallskip

We denote $\Sigma^w_{r,0}$, $\Sigma^s_{r,0}$ and $\Sigma^w_{r,\eps}$, 
$\Sigma^s_{r,\eps}$ the $r$-sewing and $r$-sliding regions of systems 
\eqref{eq1} and \eqref{eq22} respectively.
The $r$-regularizations of \eqref{eq1} and \eqref{eq22} 
are 
\begin{equation*}
\dot{x}=1/2\big(\left(1+\psi\left(h/\delta,p\right)\right)F
+\left(1-\psi\left(h/\delta,p\right)\right)G\big), \quad
\e\dot{y}=H
\end{equation*}
and
\begin{equation*}
\dot{x}=1/2\big((1+\psi(\tilde{h}/\delta,p))\tilde{F}
+(1-\psi(\tilde{h}/\delta,p))\tilde{G}\big), \quad
\tilde{H}=0,
\end{equation*} 
respectively, where $\psi:\R \times\Sigma\rightarrow [-1,1]$ is a more general smooth 
transition function satisfying that $\psi(t,p)=-1$ for $t\leq-1$ and 
$\psi(t,1)=1$ for $t\geq1$.

The next theorem states that  $r$-sliding points $p_0\in\Sigma_{r,0}^s$ 
persist under effect of singular perturbations with the additional 
assumptions:
\begin{equation}\label{H1}
\frac{\partial H}{\partial y}(p_0,0)\neq0,\quad 
\dfrac{\partial h}{\partial x} (p_0,0)\neq0, \quad 
\dfrac{\partial h}{\partial y}\equiv 0, 
\end{equation}
\begin{equation}\label{H2}
\widetilde{F}.h(p_0,0)\neq\widetilde{G}.h(p_0,0),\quad 
\frac{\partial \psi}{\partial t} (t,p_0)\neq0\quad \hbox{for}\, -1<t<1.
\end{equation}

The assumptions are necessary for applying the change of coordinates 
described in the proof and for ensuring the normal 
hyperbolicity of $p_0$.
In \cite{Ko} Sieber and Kowalczyk show that stable periodic motion with 
sliding is not robust under effect of singular perturbations.  Fridman (\cite{F1,F2}) 
also studies periodic motion considering the last assumption of \eqref{H1}.
In \cite{Ped1} the authors provide examples showing that sliding regions 
are not persistent with respect singular perturbations if this 
assumption is not considered.

\begin{theorem} \label{teoB} Consider a non--smooth slow--fast system \eqref{eq1} and 
	$p_0\in\Sigma_{r,0}^s$ satisfying the assumptions  \eqref{H1} and 
	\eqref{H2}. Then the following statements hold.
\begin{itemize}
	\item[(a)] There exist sufficiently small $\e_0>0$  and a
		family of $r$-sliding points $\{p_{\e}:\,\e\in(0,\e_0) \}$
		of system \eqref{eq1} such that $p_{\e}\rightarrow p_0$ as
		$\e\rightarrow 0$, according to Hausdorff distance.
	
	\smallskip
	
	\item[(b)] If $p_0$ is an equilibrium point (or periodic orbit) of 
		the sliding vector field associated to reduced system \eqref{eq22} then 
		there exist sufficiently small $\e_1>0$  and a family of equilibrium
		points (or periodic orbits) $\{p_{\e}:\,\e\in(0,\e_1) \}$ of the 
		sliding vector field associated to system \eqref{eq1} such that
		$p_{\e}\rightarrow p_0$ as $\e\rightarrow 0$, according to 
		Hausdorff distance.
\end{itemize}
\end{theorem}

\begin{example}\rm
	Consider the non--smooth slow--fast system
	\begin{equation}\label{ex2}
	(\dot{x}_1,\dot{x}_2)=\left\{
	\begin{array}{lcc}
	(x_2-1,-1+\e), &\hbox{if}&  x_1\geq 0,\\
	(x_1+x_2+\e,x_2+1-\e), &\hbox{if}&  x_1\leq 0,
	\end{array}
	\right.  \quad \e\dot{y}=y.
	\end{equation}
	The corresponding reduced system
	\begin{equation}\label{ex2red}
	(\dot{x}_1,\dot{x}_2)=\left\{
	\begin{array}{lcc}
	(x_2-1,-1), &\hbox{if}&  x_1\geq 0,\\
	(x_1+x_2,1), &\hbox{if}&  x_1\leq 0,
	\end{array}
	\right.
	\end{equation}
	is defined on the plane $\{y=0\}$. Note that $(0,0,0)$ and $(0,1,0)$ 
	are fold points and $\Sigma^s=\left] 0,1\right[ $. See Figure \ref{fig2}--(A). 
	
	\smallskip
	
	The $r$-regularization of \eqref{ex2} is  
	\begin{eqnarray*}
		\dot{x}_1&=&1/2\big(-1 +x_1+2 x_2+\e -(x_1+\e 
		+1)\psi({x_1}/{\delta},x_2,y)  \big),\\
		\dot{x}_2&=&1/2 \big(x_2-(2+x_2
		-2 \e )\psi({x_1}/{\delta},x_2,y)  \big),\\
		\e\dot{y}&=&y.
	\end{eqnarray*}
	Assume that the partial derivative 
	$\frac{\partial \psi}{\partial t}(t,x_2,0)$
	vanishes only at $(t,x_2,0)=(a_0,b_0,0)$ with $a_0=0$ and $b_0=-2$.
	Applying the directional blow--up $x_1=\delta \overline{x}_1$ the 
	previous system becomes
	\begin{eqnarray}
	\delta\dot{\overline{x}}_1&=&1/2\big(-1 +\delta\,\overline{x}_1+2 x_2
	+\e -(\delta\,\overline{x}_1+\e 
	+1)\psi(\overline{x}_1,x_2,y)  \big), \nonumber \\
	\dot{x}_2&=&1/2  \big(x_2-(2+x_2
	-2 \e )\psi(\overline{x}_1,x_2,y)  \big),\label{regex2}\\
	\e\dot{y}&=&y. \nonumber
	\end{eqnarray}
	The reduced system ($\e=\delta=0$) associated to \eqref{regex2} is
	\begin{eqnarray*}
		0&=&-1 +2 x_2 -\psi\left(\overline{x}_1,x_2,0\right) ,\\
		\dot{x}_2&=&1/2   \big(x_2-(2+x_2) \psi(\overline{x}_1,x_2,0)  
		\big),\\
		0&=&y,
	\end{eqnarray*}
	and it is a $r$--regularization of \eqref{ex2red}. The slow manifold 
	$\mathcal{S}=\{(\overline{x}_1,(1+\psi(\overline{x}_1,x_2,0))/2,0)\}$ 
	is a curve connecting the points $(-1,0,0)$ and $(1,1,0)$. All points 
	in $\mathcal{S}$ are normally hyperbolic (parameter $\delta$), except 
	the point $(0,-2,0)$.
	According to Fenichel's result $\mathcal{S}_0=\mathcal{S}-\{(0,-2,0)\}$  
	persists for the system
	\begin{eqnarray*}
		\delta\dot{\overline{x}}_1&=&1/2\big(-1 +\delta\,\overline{x}_1+2 x_2
		-(\delta\,\overline{x}_1
		+1)\psi(\overline{x}_1,x_2,y)  \big),  \\
		\dot{x}_2&=&\frac{1}{2}   \big(x_2-(2
		+x_2) \psi\left(\overline{x}_1,x_2,y\right)  \big),\\
		0&=&y,
	\end{eqnarray*}
	i.e, there exists an invariant manifold $\mathcal{S}_{0,\delta}$ of the 
	previous system converging to $\mathcal{S}_0$, as $\delta \rightarrow0$. 
	So, $\Sigma^s_{r,0}=\left] -2,1\right[ $.
	For each small $\delta>0$, $\mathcal{S}_{0,\delta}$ is normally hyperbolic 
	(parameter $\eps$) for the previous system. Thus, there exists an invariant 
	manifold $\mathcal{S}_{\e,\delta}$ of \eqref{regex2} such that 
	$$\mathcal{S}_{\e,\delta}\stackrel{\scriptscriptstyle{\eps\rightarrow 0}}
	{\longrightarrow} \mathcal{S}_{0,\delta}\stackrel{\scriptscriptstyle
		{\delta\rightarrow 0}}{\longrightarrow} \mathcal{S}_0,$$
	that is, there exists a $r$-sliding region $\Sigma^s_{r,\e}$ of 
	\eqref{ex2} converging to $\Sigma^s_{r,0}$.
	
	
	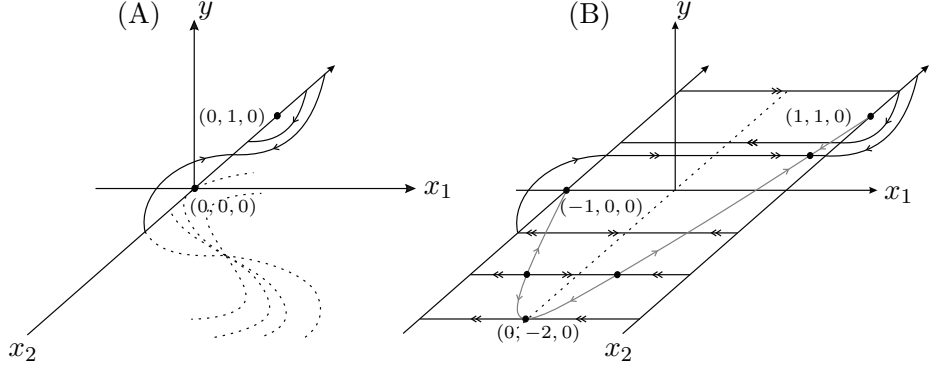
\begin{figure}
		\begin{overpic}[width=12cm]{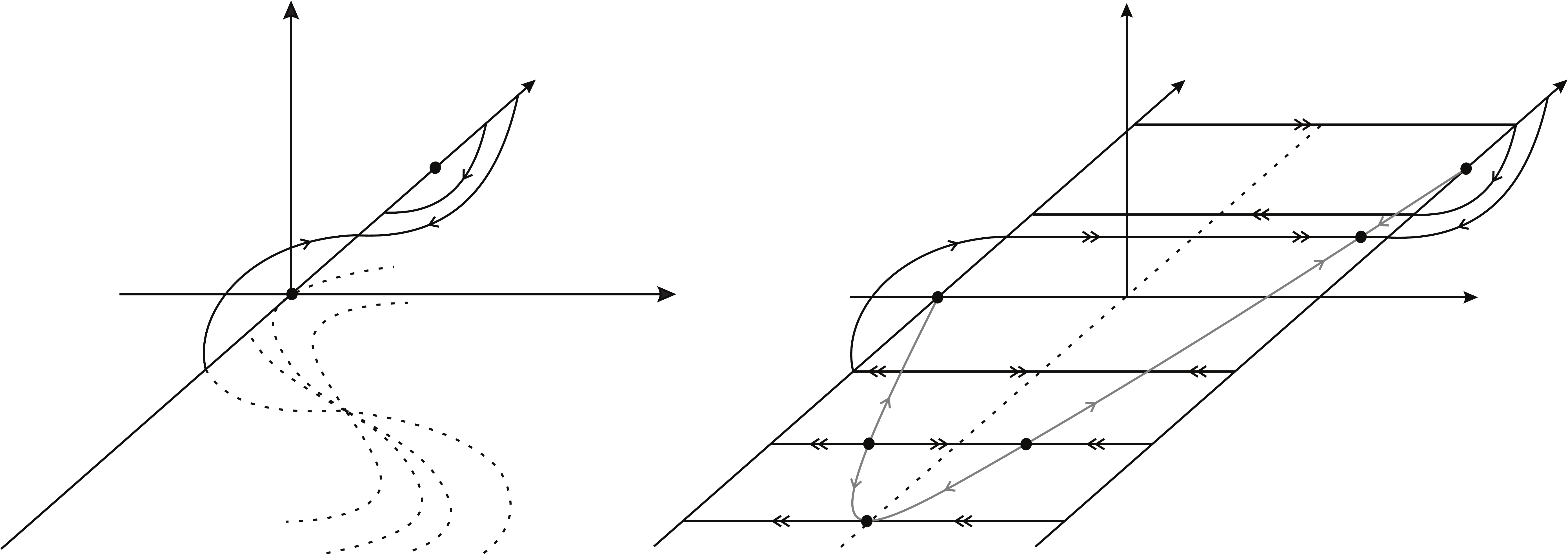}
			\put(10,35){(A)}
			\put(60,35){(B)}
			\put(44,16){$x_1$}
			\put(95,15.5){$x_1$}
			\put(-2,-2){$x_2$}
			\put(64,-2){$x_2$}
			\put(19,36){$y$}
			\put(72,36){$y$}
			\put(52,0){\tiny $ (0,-2,0)$}
			\put(18,14){\tiny $(0,0,0)$}
			\put(19,24){\tiny $(0,1,0)$}
			\put(59,14){\tiny $ (-1,0,0)$}
			\put(84,24){\tiny $(1,1,0)$}
		\end{overpic}
		\caption{\footnotesize In  (A) we have the folds points $(0,0,0)$ and 
			$(0,1,0)$. The sliding region is $\Sigma^s=\left] 0,1\right[ $. In  (B) we obtain 
			after the blow--up a smooth curve (critical manifold) connecting the 
			fold points. The point $(0,-2,0)$ is a non normal hyperbolic point.}
		\label{fig2}
	\end{figure}
	
	\smallskip
	
	The sliding vector field of \eqref{ex2red} is 
	\begin{equation}\label{slred}
	x_1=0,\quad\dot{x}_2=1-x_2-x_2^2,\quad y=0.
	\end{equation}
	It has the equilibrium points
	\[p_0^{\pm}=\left(0,{(-1\pm\sqrt{5}})/{2},0\right).\]
	The sliding vector field of \eqref{ex2} is the 
	slow--fast system
	\begin{equation*}
	x_1=0,\quad\dot{x}_2=\frac{-x_2^2+2 \e x_2 
		-x_2+\e ^2-2 \e +1}{\e
		+1},\quad \e\dot{y}=y.
	\end{equation*}
	Note that its reduced system coincides with \eqref{slred}. 
	It has the equilibrium points
	\[
	p_{\e}^{\pm}=\bigg(0,\left({-1+2\e\pm\sqrt{5-12\e+8\eps^2}}\right)/{2},0\bigg),
	\]
	converging respectively to $p_0^+$ and $p_0^-$ as $\e$ goes to zero.
\end{example}
\begin{proposition}\label{prop1}
	Consider the non--smooth slow--fast system \eqref{eq1} and  
	$p_0\in\Sigma^s_{r,0}$ satisfying the last two assumptions of \eqref{H1}.
	Then there exist local coordinates around $(p,\e)=(p_0,0)$ such 
	that $h(x,y,\e)=x_1$.
\end{proposition}	
\noindent \emph{Proof.}	
		Without lost of generality we can suppose that 
		$\partial h/\partial x_1(p_0,0)\neq0$.
		Applying the change of coordinates $\overline{x}_1=h(x,y,\e)$,
		$\overline{x}_i=x_i$ and $\overline{y}=y$, for $i=2,...,n$,
		we obtain
		\begin{eqnarray*}
			\dot{\overline{x}}_1&=&\dfrac{\partial h}{\partial {x_1}} \dot{x}_1+...+ 
			\dfrac{\partial h}{\partial {x_n}}\dot{x}_n
			+\dfrac{\partial h}{\partial {y}}\,\dot{y}+\dfrac{\partial h}
			{\partial {\e}}\, \dot{\eps}\\
			&=&\dfrac{\partial h}{\partial {x_1}} \dot{x}_1+...+ \dfrac{\partial h}
			{\partial {x_n}}\dot{x}_n+\dfrac{\partial h}{\partial {y}} \dfrac{H}{\e}
			+\dfrac{\partial h}{\partial {\e}}\, \dot{\eps}.
		\end{eqnarray*}
		Since $\partial h/\partial y\equiv0$ and $\dot{\eps}=0$, the previous 
		expression becomes 
		\[\dot{\overline{x}}_1=\dfrac{\partial h}{\partial {x_1}} \dot{x}_1+...
		+ \dfrac{\partial h}{\partial {x_n}}\dot{x}_n.\]
		Since $\partial h/\partial {x_1}\neq0$ the determinant of the change of 
		coordinates matrix is nonzero and system \eqref{eq1} becomes
		\begin{equation*}
		\dot{\overline{x}}=\left\{
		\begin{array}{ll}
		\overline{F}(\overline{x},\overline{y},\e), \quad\hbox{if} \quad \overline{x}_1\geq 0\\
		\overline{G}(\overline{x},\overline{y},\e),\quad \hbox{if} \quad \overline{x}_1\leq 0
		\end{array}
		\right. , \quad \e\dot{\overline{y}}=H(\overline{x},\overline{y},\e).
		\end{equation*}
		where $\overline{x}=(\overline{x}_1,...,\overline{x}_n)$.\bbox\\

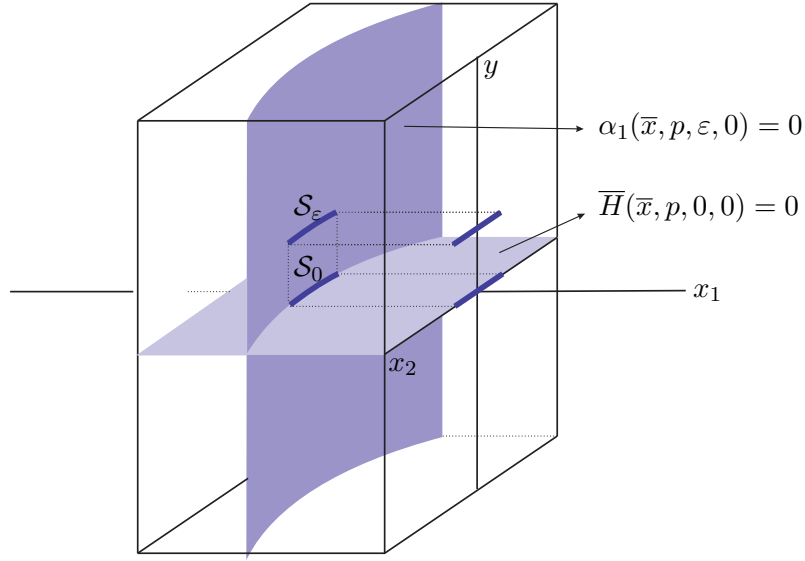
\begin{figure}
	\begin{overpic}[width=9cm]{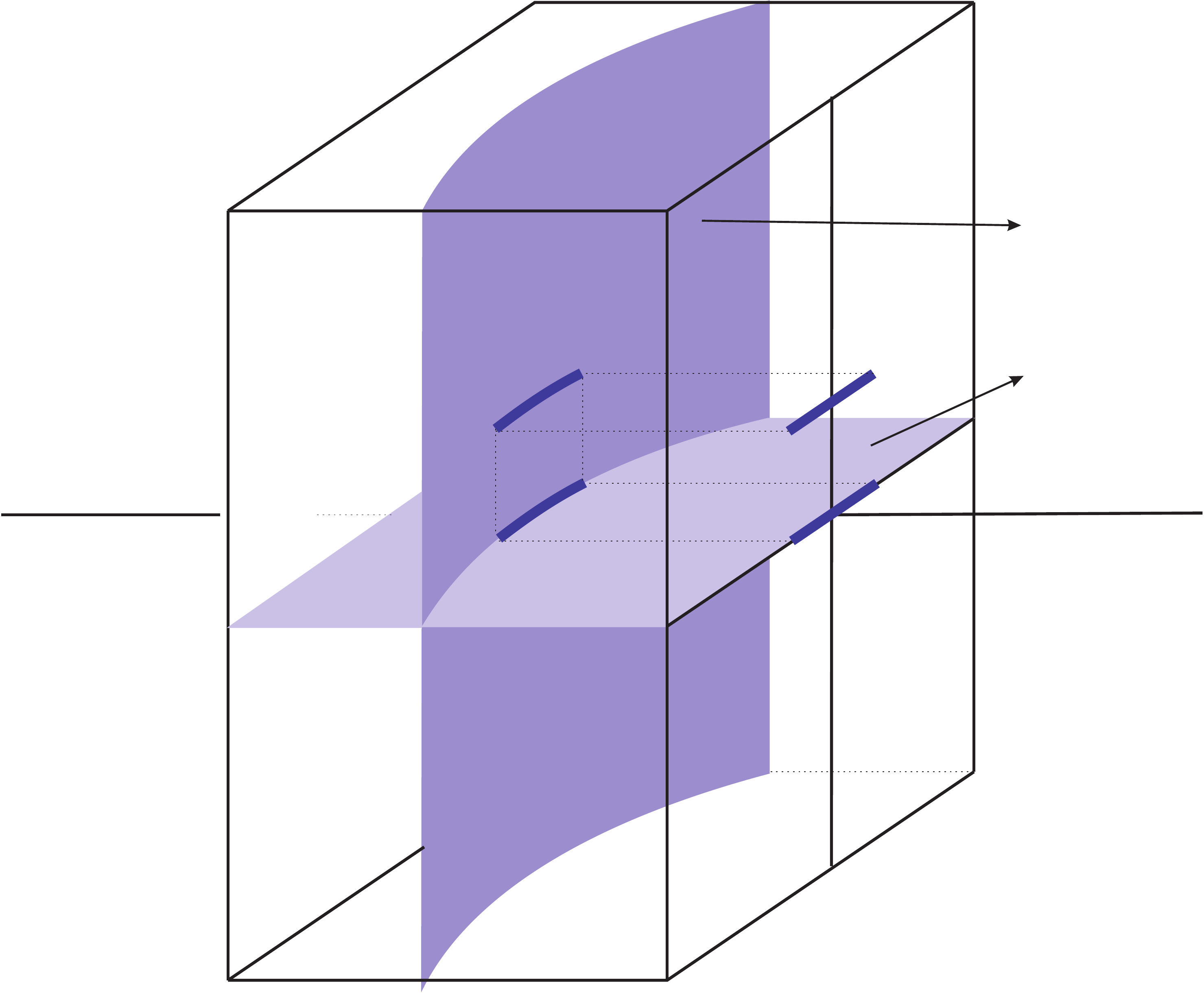}
		\put(42,51){$\mathcal{S}_{\eps}$}
		\put(42,42){$\mathcal{S}_{0}$}
		\put(101,39){$x_1$}
		\put(70,72){$y$}
		\put(56,28){$x_2$}
		\put(87,63){$\alpha_1(\overline{x},p,\eps,0)=0$}
		\put(87,51){$\overline{H}(\overline{x},p,0,0)=0$}
	\end{overpic}
	\caption{\footnotesize Geometric situation in the blow--up process 
		described in the proof of Theorem \ref{teoB}.}
	\label{figteoB}
\end{figure}

\noindent{\textit{Proof of Theorem \ref{teoB}}}. Suppose that 
$p_0\in\Sigma^s_{r,0}$. According to Proposition \ref{prop1} we can 
assume that $h(x,y,\e)=x_1$ in \eqref{eq1} around $(p,\e)=(p_0,0)$. 
Denote $p=(x_2,...,x_n,y)$. The switching manifold becomes 
$\Sigma_\eps=\{(0,p)\}$ for each  $\eps\geq 0$. The $r$-regularization 
of \eqref{eq1} is the $2-$parameters ($\eps$ and $\delta$) family
\begin{equation*}\label{regB}
\dot{x}=1/2\big((1+\psi({x_1}/{\delta},p))F
+(1-\psi({x_1}/{\delta},p))G\big),\quad
\e\dot{y}=H.
\end{equation*}
Note that it is a slow--fast system (parameter $\e$) and its
reduced system is a $r$-regularization of reduced system
\eqref{eq22}. Applying the directional blow--up 
$x_1=\delta \overline{x}_1$ we obtain the \emph{three time
	scale singular perturbation problem}
\begin{eqnarray}
\delta\dot{\overline{x}}_1&=&\alpha_1(\overline{x}_1,p,\e,\delta),\nonumber\\
\dot{x}_i&=&\alpha_i(\overline{x}_1,p,\e,\delta),\label{eqth1}\\
\e\dot{y}&=&\overline{H}(\overline{x}_1,p,\e,\delta),\nonumber
\end{eqnarray}
for $i=2,...,n$, where
\[\alpha_i(\overline{x}_1,p,\e,\delta)=1/2\big((1
+\psi(\overline{x}_1,p))F_i+(1-\psi(\overline{x}_1,p))G_i\big),\]
for $i=1,...,n$, with the function $F_i$ and $G_i$ evaluated at 
$(\delta \overline{x}_1,p,\eps)$ and $\overline{H}(\overline{x}_1,p,\e,\delta)
=H(\delta \,\overline{x}_1, p,\e)$.
The reduced system associated to \eqref{eqth1} ($\e=\delta=0$) 
is the following 
\begin{eqnarray}
0&=&\alpha_1(\overline{x}_1,p,0,0),\label{eqcri1}\\
\dot{x}_i&=&\alpha_i(\overline{x}_1,p,0,0),\nonumber \\
0&=&\overline{H}(\overline{x}_1, p,0,0),\label{eqcri2}
\end{eqnarray}
for $i=2,...,n$. Let $\mathcal{S}$ be the critical manifold given by 
equations \eqref{eqcri1} and \eqref{eqcri2} and 
$(t_0,p_0)\in\mathcal{S}$. Since
\[\dfrac{\partial \overline{H}}{\partial y}(t_0, p_0,0,0)
=\dfrac{\partial H}{\partial y}(0,p_0,0)\neq0,\]
Proposition \ref{Fenichel} says that the point $(t_0,p_0)$
persists for the system
\begin{eqnarray}
0&=&\alpha_1(\overline{x}_1,p,\e,0),\nonumber\\
\dot{x}_i&=&\alpha_i(\overline{x}_1,p,\e,0),\label{eqth}\\
\e\dot{y}&=&\overline{H}(\overline{x}_1,p,\e,0),\nonumber
\end{eqnarray}
for $i=2,...,n$. Indeed there exists a compact set 
$\mathcal{S}_0\subset\mathcal{S}$ containing $(t_0,p_0)$ and a 
family of invariant manifolds $\mathcal{S}_{\e}$ of \eqref{eqth},
such that $\mathcal{S}_{\e}\rightarrow\mathcal{S}_0$ as 
$\e \rightarrow 0$, according to Hausdorff distance. See Figure 
\ref{figteoB}. Since
\[\dfrac{\partial \alpha_1}{\partial \overline{x}_1}(\overline{x}_1, p,0,0)=
\dfrac{\partial \psi}{\partial \overline{x}_1}(\overline{x}_1,p).(F_1-G_1)(0,p,0)\neq0,\]
for $(\overline{x}_1,p)\in\mathcal{S}_0$, by continuity we have that 
$$\dfrac{\partial \alpha_1}{\partial \overline{x}_1}(\overline{x}_1, p,0,\e)\neq0,$$
for $(\overline{x}_1,p)\in\mathcal{S}_{\e}$ and sufficiently small 
$\eps>0$. Thus Proposition \ref{Fenichel} ensures the persistence of
$\mathcal{S}_{\e}$ for system \eqref{eqth1}. More specifically, there 
exists a family of invariant manifolds $\mathcal{S}_{\e,\delta}$ of 
\eqref{eqth1} such that 
$\mathcal{S}_{\e,\delta}\rightarrow \mathcal{S}_{\e}$ as
$\delta\rightarrow 0$.

\smallskip

Therefore there exists a family of $r$-sliding points $p_{\e}$ 
of \eqref{eq1} satisfying that $p_{\e}\rightarrow p_0$
as $\e\rightarrow 0$, according to Hausdorff's distance. It concludes
the proof of statement (a).

\smallskip

Note that the sliding vector field associated to \eqref{eq1}
is the slow--fast system
\begin{equation*}
\dot{x}_i =\dfrac{F_1G_i - G_1F_i}{F_1-G_1}(0,p,\e),\quad \e\dot{y}=H(0,p,\e),
\end{equation*}
for $i=2,...,n$, and its reduced system 
\begin{equation*}
\dot{x}_i =\dfrac{F_1G_i - G_1F_i}{F_1-G_1}(0,p,0),\quad H(0,p,0)=0,
\end{equation*}
has the dynamics of the sliding vector field associated to \eqref{eq22}.
So the proof of statement (b) follows directly from Proposition \ref{Fenichel}.\bbox

\section{Continuous Combinations of Non--Smooth Systems}\label{s5}

Now we consider another way to define sliding points. Instead of
considering a convex combination of vectors $X^+(p)$ and $X^-(p)$  we consider a 
\textit{continuous combination}. This convention is given in 
\cite{NJ} and \cite{DIP}.

\smallskip

Consider a non--smooth system \eqref{sis1}. A continuous combination 
of $X^+$ and $X^-$ is a 1--parameter family of smooth vector fields 
$\widetilde{X}(\lambda,p)$ with $(\lambda,p) \in [-1,1]\times \Sigma$  
satisfying that $\widetilde{X}(1,p)=X^+(p)$ and 
$\widetilde{X}(-1,p)=X^-(p)$. We denote
\[ [X^+,X^-]^c=\{\widetilde{X}(\lambda,p), \lambda \in [-1,1]\}.\]

Consider and a regular point $p\in\Sigma$.
\begin{itemize}
\item[(i)] We say that  $p$ is a\emph{ $c$-sewing} point and denote
$p\in\Sigma^{w}_{c}$ if  $\widetilde{X}.h(\lambda,p)\neq 0 $ for all
$\lambda\in (-1,1).$

\smallskip

\item[(ii)] We say that  $p$ is a \emph{$c$-sliding} point and denote
$p\in\Sigma^{s}_c$ if there exists  $\lambda\in (-1,1)$, such that 
$\widetilde{X}.h(\lambda,p)=0$.
\end{itemize}
We say that $\widetilde{X}(\lambda (p),p)$ is
a \emph{$c$-sliding vector field} if for each $p\in \Sigma^{s}_c$
there exists $\lambda (p)\in (-1,1)$ such that  
$\widetilde{X}.h(\lambda(p),p)=0$.

\smallskip

\begin{figure}[h]
\begin{overpic}[width=6cm]{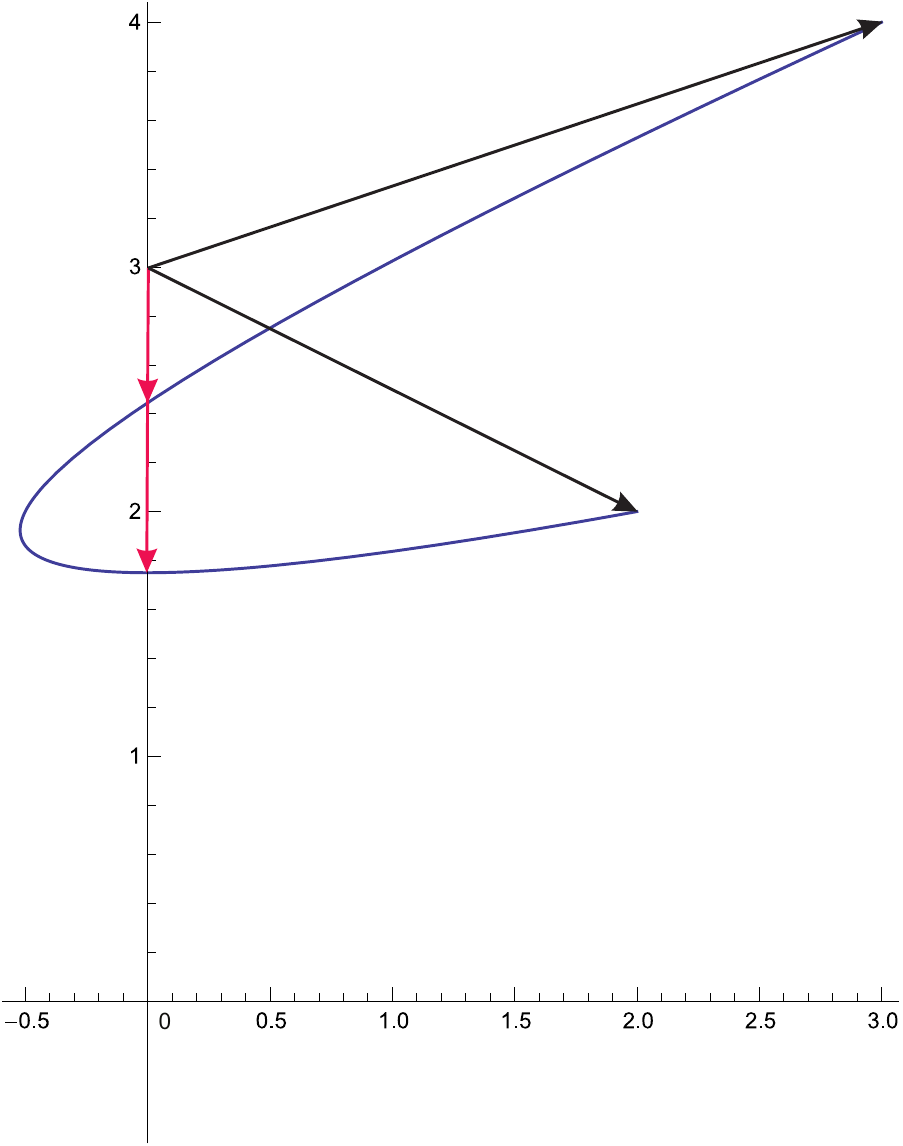}
\put(80,12){$x_1$}
\put(10,102){$x_2$}
\put(38,92){\small $ X^{-}(0,3)$}
\put(27,58){\small $X^{+}(0,3)$}
\put(-8,53){\small $X_{0}^{\lambda_1}(0,3)$}
\put(-8,70){\small $X_{0}^{\lambda_2}(0,3)$}
    \end{overpic}
\caption{\footnotesize The curve in the continuous
combinations of the vectors (in red) $X^{-}(0,3)$ and $X^{+}(0,3)$ 
in Example \ref{ex1}. The vectors in red represent the $c-$sliding 
vector fields $X_{0}^{\lambda_1}$ and $X_{0}^{\lambda_2}$
at $(0,3)$.}
\label{fig1}
 \end{figure}
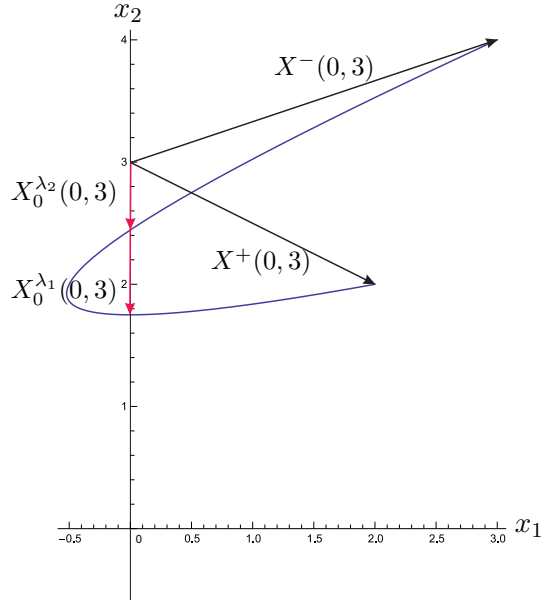

There may be more than one possible sliding on $p$, see Figure \ref{fig1}. 
In \cite{DIP} 
the authors prove that $\Sigma_c^w\subseteq\Sigma^w$ and 
$\Sigma^s\subseteq\Sigma^s_c$.

\smallskip

Here we study the persistence of $c$-sliding
points via slow--fast systems considering the continuous combination 
\[\widetilde{Y}(\lambda,x,y,\e)=
(\widetilde{X}(\lambda,x,y,\e),H(x,y,\e)/\e),\]
where
$\widetilde{X}(\lambda,x,y,\e)$ is a continuous combination
of $X^{+}(x,y,\e)$ and $X^{-}(x,y,\e)$. 

\smallskip

We denote $\Sigma^w_{c,0}$, $\Sigma^s_{c,0}$ and $\Sigma^w_{c,\eps}$, 
$\Sigma^s_{c,\eps}$ the $c$-sewing and $c$-sliding regions of systems 
\eqref{eq1} and \eqref{eq22} respectively.

\smallskip

The next theorem provides results like the ones given in Theorem \ref{teoB}
however for $c$-sliding points. We consider the assumption \eqref{H1} 
and the following one
\begin{equation}\label{H3}
\frac{\partial Q}{\partial \lambda} (\lambda^*,p_0,0)\neq0,
\end{equation}
where $Q(\lambda,p,\eps)=\widetilde{Y}.\nabla h(\lambda,p,\eps)$ and 
$\lambda^*$ satisfies the equation $Q(\lambda^*,p_0,0)=0$ for 
$p_0\in\Sigma^s_{c,0}$.

\begin{theorem}\label{teoC} Consider a non--smooth slow--fast system \eqref{eq1} and 
$p_0\in\Sigma_{r,0}^s$ satisfying the assumptions  \eqref{H1} and 
\eqref{H3}.Then the following statements hold.
\begin{itemize}
\item[(a)]  There exist sufficiently small $\e_0>0$  and a
family of $c$-sliding points $\{p_{\e}:\,\e\in(0,\e_0) \}$
of system \eqref{eq1} such that $p_{\e}\rightarrow p_0$ as
$\e\rightarrow 0$, according to Hausdorff distance.
	
\smallskip
	
\item[(b)] If $p_0$ is an equilibrium point (or periodic orbit) 
of the $c$-sliding vector field associated to reduced system 
\eqref{eq22} then there exist sufficiently small $\e_1>0$  and a 
family of equilibrium points (or periodic orbits) 
$\{p_{\e}:\,\e\in(0,\e_1) \}$ of the $c$-sliding vector field 
associated to system \eqref{eq1} such that
$p_{\e}\rightarrow p_0$ as $\e\rightarrow 0$,
according to Hausdorff distance.
\end{itemize}
\end{theorem}

\begin{example}\label{ex1}\rm
Consider the non--smooth slow--fast system
\begin{equation}\label{sisex33}
\dot{x}=\left\{
  \begin{array}{ccc}
    (x_2-1+\e,-1+\e), &\hbox{if}&  x_1\geq 0\\
    (x_2+\e,1+\e), &\hbox{if}&  x_1\leq 0
  \end{array}
\right. , \quad \e\dot{y}=y
\end{equation}
and the continuous combination
\begin{equation*}\label{eqex1}
\widetilde{Y}(\lambda,x_1,x_2,y,\e)=\big(\lambda^2(x_2+\e)
-(\lambda+1)/2,\lambda^2-\lambda-1+\e,y/\e\big).
\end{equation*}
The equation $\lambda^2(x_2+\e)-(\lambda+1)/2=0$ provides 
\[\lambda_{1}^\eps=\dfrac{1+\sqrt{1+8x_2+8\e}}{4(x_2+\e)},
\quad \lambda_{2}^\eps=\dfrac{1+\sqrt{1+8x_2+8\e}}{4(x_2+\e)} .\]
Replacing $\lambda_{1}^\eps$ and $\lambda_{2}^\eps$ in 
$\widetilde{Y}$ we obtain two $c$-sliding vector fields
\[X^{\lambda_1^\eps}_{\e}=\left(0,\frac{\left(\sqrt{8 x_2+8 \e
+1}+1\right)^2}{16 (x_2+\e )^2}-\frac{\sqrt{8 x_2+8 \e 
+1}+1}{4 (x_2+\e )}+\e -1,\frac{y}{\e }\right)\]
and
\[X^{\lambda_2^\eps}_{\e}=\left(0,\frac{\left(\sqrt{8 x_2+8 \e
+1}-1\right)^2}{16 (x_2+\e)^2}+\frac{\sqrt{8 x_2+8 \e
+1}-1}{4 (x_2+\e )}+\e -1,\frac{y}{\e }\right).\]

Note that $x_2\geq (-1-8\e)/8$. For $|\lambda_1|<1$ and
$|\lambda_2|<1$ we obtain
$\Sigma^{s}_{c,\e}=\left] -\e,1-\e\right[ \cup   \left] 1-\e,+\infty \right[  \times \R$.
In $\left]1-\e,+\infty\right[\times \R$  are defined two $c$--sliding vector
fields ($X^{\lambda_1}_\e$ and $X^{\lambda_2}_\e$) and in 
$\left]-\e,1-\e\right[\times\R$ it is only defined $X^{\lambda_2}_\e$.
The dynamics of $X_{\e}^{\lambda_1}$ and $X_{\e}^{\lambda_2}$ 
are governed respectively by slow--fast systems 
\[(\dot{x}_1,\dot{x}_2,\dot y )=X^{\lambda_1^\eps}_{\e},
\quad(\dot{x}_1,\dot{x}_2,\dot y )=X^{\lambda_2^\eps}_{\e}. \]
The reduced systems are given respectively by  
\begin{equation}\label{reduz1}
(\dot{x}_1,\dot{x}_2)=\left(0,-\frac{(2 x_2-1) \left(4 x_2+\sqrt{8
x_2+1}+1\right)}{8 x_2^2}\right),\quad y=0,
\end{equation}
and
\begin{equation}\label{reduz2}
(\dot{x}_1,\dot{x}_2 )=\left(0,\frac{(2 x_2-1) \left(-4 x_2+\sqrt{8
x_2+1}-1\right)}{8 x_2^2}\right), \quad y=0.
\end{equation}

Now, consider the reduced system of \eqref{sisex33}
\begin{equation*}
\dot{x}=\left\{
\begin{array}{lcc}
(x_2-1,-1), &\hbox{if}&  x_1\geq 0\\
(x_2,1), &\hbox{if}&  x_1\leq 0
\end{array}
\right. ,
\end{equation*}
defined in the plane $y=0$ 
and the continuous combination
\begin{equation*}
\widetilde{X}(\lambda,x,y)=\big(\lambda^2x_2-(\lambda+1)/2,\lambda^2-
\lambda-1).
\end{equation*}
Using the same process to get $\Sigma^{s}_{c,\e}$ we obtain that
\[\lambda_{1,2}=\dfrac{1\pm\sqrt{1+8x_2}}{4x_2},\]
and $\Sigma^{s}_{c,0}=\left]0,1\right[\cup\left]1,+\infty\right[$.
In $\left]1,+\infty\right[$ are defined two $c$--sliding vector fields 
$X^{\lambda_1}_0$ and $X^{\lambda_2}_0$ and they have the dynamics 
described respectively by \eqref{reduz1} and \eqref{reduz2}. 
In $\left] 0,1\right[ $ it is only defined
$X^{\lambda_2}_0$.
The sliding vector fields $X_0^{\lambda_1}$ and $X_0^{\lambda_2}$
have a normally hyperbolic equilibrium point $p_0=1/2$. 
The $c-$sliding vector fields $X_{\e}^{\lambda_1}$ and 
$X_{\e}^{\lambda_2}$ have the same equilibrium points
\[p^{\pm}_{\e}=\frac{-4 \e ^3+8 \e ^2-5 \e +2\pm\sqrt{5 \e ^2-4 \e^3}}
{4 \left(\e ^2-2 \e +1\right)},\]
satisfying that $p^{\pm}_{\e} \rightarrow p_0$ as $\e \rightarrow 0$.
\end{example}

\noindent\textit{Proof of Theorem \ref{teoC}.} Consider the non--smooth 
system \eqref{eq1} and $p_0\in\Sigma_{c,0}^s$. 
According to Proposition \ref{prop1} 
we can assume that $h(x,y,\e)=x_1$ 
around $(p,\e)=(p_0,0)$. Denote $p=(x_2,...,x_n,y)$ and  
\[\widetilde{X}(\lambda, x_1,p)=(s_1(\lambda, x_1, p),...,
s_n(\lambda, x_1, p))\] a continuous combinations of $\widetilde{F}$
and $\widetilde{G}$.

\smallskip

The $c$-sliding regions of \eqref{eq1} and \eqref{eq22} become 
\begin{eqnarray*}
\Sigma_{c,\e}^{s}&=&\{(0,x_2,...,x_n,y): \exists \lambda\in(-1,1),\,
s_1(\lambda, x, y,\eps)=0\},\\
\Sigma_{c,0}^{s}&=&\{(0,x_2,...,x_n,y):\exists {\lambda^*}\in(-1,1),\,
s_1({\lambda^*},x,y,0)=0\}.
\end{eqnarray*}

Consider ${\lambda^*}$ such
that $s_1({\lambda^*},p_0,0)=0$. Assumption \eqref{H3} ensures
the existence of a neighborhood $V$ of $(p_0,0)$ such that
$\lambda=\lambda(q,\e)$, $\lambda(p_0,0)={\lambda^*}$ and 
$s_1(\lambda(q,\e),q,\e)=0$ for all $(q,\e)\in V$, in particular 
for $q\in V\cap\Sigma_{c,\e}^{s}$.
Therefore there exists $p_\eps\in\Sigma_{c,\e}^{s}$ such that 
$p_\eps\rightarrow p_0$ as $\eps\rightarrow 0$ and statement (a) 
is proved.

\smallskip

The dynamics of the $c$-sliding vector field $\widetilde{Y}$ is
given by  
\begin{equation}\label{eq4}
\dot{x}_i=s_i(\lambda,x,y,\e),\quad \eps\dot{y}=H(x,y,\e),\quad i=2,...,n,
\end{equation}
with $\lambda$ satisfying the equation $s_1(\lambda,x,y,\e)=0$.
Note that system \eqref{eq4} is a singular perturbation problem.
The dynamic of the $c$-sliding vector field associated to system 
\eqref{eq22} is given by 
\begin{equation}\label{eq6}
\dot{x}_i=s_i(\lambda^*,x,y(x),0),\quad i=2,...,n,
\end{equation}
with $\lambda^*$ satisfying $s_1(\lambda^*,x,y(x),0)=0$.

\smallskip

Since $\lambda=\lambda(q,\e)$ the reduced system associated to 
system \eqref{eq4} is given by
\begin{equation}\label{eq5}
\dot{x}_i=s_i(\lambda(q,0),x,y(x),0),\quad i=2,...,n,
\end{equation}
Therefore system \eqref{eq6} coincides with system \eqref{eq5}. 
Now we apply the Fenichel's result for concluding the proof of 
item (b). \bbox\\

\noindent \textbf{Remark.} In \cite{NJ} and \cite{DIP} is defined a regularization of 
\eqref{sis1} called \emph{nonlinear regularization} as 
\begin{equation}\label{nonl}
\dot{x}=\widetilde{X}\left(\varphi(h/\delta),p\right), 
\end{equation}
where $\varphi$ is a monotonic transition function. 
Next lemma says that $\Sigma_c^s$ is the 
sliding region linked to nonlinear regularization \eqref{nonl}. 
Thus Theorem \ref{teoC} can be proved such as Theorem \ref{teoB} but considering the nonlinear 
regularization of systems \eqref{eq1} and \eqref{eq22} 
\begin{eqnarray*}
	\dot{x}&=&\widetilde{Y}(\varphi(h/\e),x,y,\e_2), \\
	\dot{x}&=&\widetilde{X}(\varphi(\tilde{h}/\e),x,y(x),0),
	\quad \widetilde{H}(x)=0,
\end{eqnarray*} 
respectively.

\begin{lemma} 
Consider a non--smooth system \eqref{sis1} with 
$K(\lambda,p)=\widetilde{X}.h(\lambda,p)$ and $\widetilde{X}$ a 
continuous combination of $X^+$ and $X^-$. Suppose that there 
exists $\lambda^*\in(-1,1)$ such that 
\begin{equation}\label{eq1lemma}
\dfrac{\partial K}{\partial \lambda}(\lambda^*,p_0)\neq0,
\end{equation}
for $p_0\in\Sigma$. 
Thus $p_0$ is a $c$-sliding point if and only if $p_0$ is a 
sliding point for nonlinear regularization \eqref{nonl}.
\end{lemma}

\noindent \emph{Proof.}	
Consider a system like \eqref{sis1} and a continuous combination 
$\widetilde{X}$ of $X^{+}$ and $X^{-}$. Take local coordinates 
$(x_1,p)=(x_1,...,x_n)$ such that $h(x)=x_1$ and 
$\widetilde{X}(\lambda,x)=(s_1(\lambda,x_1,p),..., s_n(\lambda,x_1,p))$.
The switching manifold and the $c$-sliding region become
\[\Sigma=\{(0,p)\},\quad\Sigma_{c}^{s}=\{(0,p): \exists \lambda\in(-1,1),
\,s_1(\lambda, 0, p)=0\},\]
respectively. Consider $p_0\in\Sigma^s_c$. Thus there exists $\lambda^*$ 
satisfying $s_1(\lambda^*,0,p)=0$. 
The nonlinear regularization \eqref{nonl} becomes
\begin{equation}\label{nonl1}
\dot{x}_1=s_1(\varphi(x_1/\delta),x_1,p),\quad \dot{x}_i=
s_n(\varphi(h/\delta),x_1,p), 
\end{equation}
for $i=2,...,n$. Taking the blow--up $x_1=\delta\overline{x}_1$ 
the previous system becomes the slow--fast system
\[\delta\dot{\overline{x}}_1=s_1(\varphi(\overline{x}_1),
\delta\overline{x}_1,p),\quad \dot{x}_i=s_n(\varphi(\overline{x}_1),
\delta\overline{x}_1,p),\]
for $i=2,...,n$. Define $\lambda=\varphi(\overline{x}_1)$. Since 
$\varphi'(t)>0$ in $(-1,1)$ there exists $\overline{x}_1^*$ such that 
$\lambda^*=\varphi(\overline{x}_1^*)$ with 
$(\overline{x}_1^*,p_0)\in\mathcal{M}$, where $\mathcal{M}=\{s_1(\lambda,0,p)=0\}$ is the 
slow manifold. We claim that 
$(\overline{x}_1^*,p_0)$ is a normally hyperbolic point. In fact,
\[\dfrac{\partial s_1 }{\partial\, \overline{x}_1}(\lambda^*,0,p_0)=
\varphi'(\overline{x}_1^*)\, \dfrac{\partial s_1}{\partial \lambda} 
(\lambda^*,0,p_0)\neq0.\] 
Thus the Fenichel's result ensures the existence of an invariant 
manifold $\mathcal{M}_\lambda$ of \eqref{nonl1} converging to a compact 
manifold $\mathcal{M}_0\subset\mathcal{M}$  containing $p_0$. Therefore 
$p_0$ is a sliding point for the nonlinear regularization \eqref{nonl}.

\smallskip

Conversely, if $p_0$ is a sliding point for the nonlinear regularization 
\eqref{nonl} satisfying assumption \eqref{eq1lemma} then $p_0$ is a 
$c$-sliding point. In fact, note that the slow manifold 
$\mathcal{M}$ and $\Sigma_c^s$ are 
defined by same equation and $\varphi$ is increasing in $(-1,1)$. So, 
there exist neighborhoods $U\subset \mathcal{M}$ and $V\subset\Sigma_c^s$ 
of $(\lambda^*,p_0)$ and $p_0$, respectively, and a diffeomorphism 
$\xi:U\rightarrow V$ satisfying $\xi(\lambda^*,p_0)=p_0$.\bbox

\section{Acknowledgments}
Jaime R. de Moraes is partially supported by FUNDECT--219/2016. 
Paulo R. da Silva is partially supported by CAPES and FAPESP.

\end{document}